\documentclass[10pt]{article}
\usepackage{latexsym,amsmath,amssymb,graphics,amscd}

\makeatletter
\@addtoreset{figure}{section}
\def\thefigure{\thesection.\@arabic\c@figure}
\def\fps@figure{h,t}
\@addtoreset{table}{bsection}

\def\thetable{\thesection.\@arabic\c@table}
\def\fps@table{h, t}
\@addtoreset{equation}{section}

\makeatother

\begin{document}

\newtheorem{theorem}{Theorem}[section]
\newtheorem{definition}[theorem]{Definition}
\newtheorem{lemma}[theorem]{Lemma}
\newtheorem{remark}[theorem]{Remark}
\newtheorem{proposition}[theorem]{Proposition}
\newtheorem{corollary}[theorem]{Corollary}
\newtheorem{example}[theorem]{Example}
\newtheorem{examples}[theorem]{Examples}

\newcommand{\bfi}{\bfseries\itshape}

\newsavebox{\savepar}
\newenvironment{boxit}{\begin{lrbox}{\savepar}
\begin{minipage}[b]{15.8cm}}{\end{minipage}
\end{lrbox}\fbox{\usebox{\savepar}}}

\makeatletter
\title{{\bf Symmetry and symplectic reduction}}
\author{Juan-Pablo Ortega$^{1}$ and  Tudor S. Ratiu$^{2}$}
\addtocounter{footnote}{1}
\footnotetext{Centre National de la Recherche Scientifique,
D\'epartement de Math\'ematiques de Besan\c con,
Universit\'e de Franche-Comt\'e.
UFR des Sciences et Techniques.
16, route de Gray.
F-25030 Besan\c con cedex. France. {\texttt
Juan-Pablo.Ortega@math.univ-fcomte.fr}. }
\addtocounter{footnote}{1}
\footnotetext{Section de Math\'ematiques, 
\'Ecole Polytechnique F\'ed\'erale de Lausanne,  CH-1015 Lausanne,
Switzerland. {\texttt Tudor.Ratiu@epfl.ch}.}
\date{}
\makeatother
\maketitle

\begin{abstract}
This encyclopedia article briefly reviews without proofs some of the main results
in symplectic reduction. The article recalls most the necessary prerequisites to
understand the main results, namely, group actions, momentum maps, and coadjoint
orbits, among others.
\end{abstract}

\section{Introduction}

The use of symmetries in the
quantitative and qualitative study of dynamical systems has a long
history that goes back to the founders  of mechanics. In most
cases, the symmetries of a system are used to implement a
procedure  generically known under the name of {\bfi  reduction}
that restricts the study of its dynamics to a system of
smaller dimension. This procedure is also used in a purely
geometric context to construct new nontrivial manifolds having various
additional structures.

Most of the reduction methods 
can be seen as constructions that systematize the  techniques
of elimination of variables found in classical mechanics.
These procedures consist basically of two steps. First, one
restricts the dynamics to flow invariant submanifolds of the system
in question and secondly one projects the restricted dynamics onto 
the symmetry orbit quotients of the spaces
constructed in the first step. Sometimes, the flow
invariant manifolds appear as the level sets of a
{\bfi  momentum map} induced by the symmetry of the system.

\section{Symmetry reduction}
\label{Symmetry reduction}

\noindent {\bf The symmetries of a system.} The standard mathematical fashion to 
describe the symmetries of a dynamical system $X \in \mathfrak{X} (M) $ 
defined on a manifold $M$ ($\mathfrak{X}(M)$ denotes the Lie algebra of
smooth vector fields on $M$ endowed with the Jacobi-Lie bracket $[\cdot ,
\cdot ]$)  consists in studying its invariance properties with respect to
a smooth Lie group
$\Phi: G
\times  M
\rightarrow M $ ({\bfi  continuous symmetries}) or Lie algebra $\phi:
\mathfrak{g}
\rightarrow \mathfrak{X} (M)  $ ({\bfi  infinitesimal symmetry}) action.
Recall that $\Phi$ is a (left) action if the map 
$g \in G \mapsto \Phi(g, \cdot ) \in \operatorname{Diff}(M)$ is a group
homomorphism, where $\operatorname{Diff}(M) $ denotes the group of smooth
diffeomorphisms of the manifold $M $. The map $\phi$ is a (left) Lie
algebra action if the map $\xi \in \mathfrak{g} \mapsto \phi( \xi) \in
\mathfrak{X}(M)$ is a Lie algebra anti-homomorphism and the map
$(m, \xi) \in M \times \mathfrak{g} \mapsto \phi( \xi)(m) \in TM$ is
smooth. The vector field $X$ is said to be $G$-{\bfi symmetric} whenever
it is equivariant with respect to the $G$-action $\Phi$, that is, $X \circ
\Phi_g=T \Phi _g \circ X $, for any $g \in  G $. The space of
$G$-symmetric vector fields on $M$ is denoted by $\mathfrak{X} (M)^G $. 
The flow $F _t $ of a $G$-symmetric vector field $X \in \mathfrak{X}
(M)^G$ is $G$-equivariant, that is, $F _t \circ  \Phi_g = \Phi_g \circ
F_t$, for any $g \in  G $. The vector field $X $ is said to be
$\mathfrak{g} $-{\bfi  symmetric} if
$[\phi(\xi), X]=0$, for any $\xi\in \mathfrak{g}$. 

If $\mathfrak{g}$ is the Lie algebra of the Lie group $G$ then the
{\bfi  infinitesimal generators} $\xi_M \in  \mathfrak{X} (M)$ of a 
smooth  $G$-group action defined by $\xi_M (m):=
\left.\frac{d}{dt}\right|_{t=0} \Phi(\exp t \xi,m)$,
$\xi\in  \mathfrak{g} $, $m \in  M $, constitute a smooth Lie algebra
$\mathfrak{g}$-action and we denote in this case $\phi( \xi) = \xi_M$.

If $m \in M$, the closed Lie subgroup $G_m: = \{ g \in G \mid \Phi(g, m)
= m\} $ is called the {\bfi isotropy\/} or {\bfi symmetry\/} subgroup of
$m$. Similarly, the Lie subalgebra $\mathfrak{g}_m : = \{ \xi\in
\mathfrak{g} \mid \phi( \xi)(m) = 0 \}$ is called the {\bfi isotropy\/}
or  {\bfi symmetry\/} subalgebra of $m$. If $\mathfrak{g}$ is the Lie
algebra of $G$ and the Lie algebra action is given by the
infinitesimal generators, then $\mathfrak{g}_m $ is the Lie algebra of
$G_m$. The action is called {\bfi free} if $G_m = \{e\} $ for every $m
\in M $ and {\bfi locally free} if $\mathfrak{g}_m = \{0\} $ for every $m
\in M $. We will write interchangeably $\Phi(g, m) = \Phi_g(m) =
\Phi^m(g) = g \cdot m $, for $m \in M $ and $g \in G $.

In this article we will focus mainly on continuous symmetries induced by 
proper Lie group actions. The action $\Phi$ is called {\bfi  proper}  
whenever for any two convergent sequences
$\{ m_{n}\}_{n \in \mathbb{N}}$ and
$\{ g_{n}\cdot m_{n}:= \Phi(g _n, m _n)\}_{n \in \mathbb{N}}$ in $M$,
there exists a convergent subsequence $\{ g_{n_{k}}\}_{k \in \mathbb{N}}$
in $G$.  Compact group actions are obviously proper.

\medskip

\noindent \textbf{Symmetry reduction of vector fields.}
Let $M$ be a smooth manifold and $G$ a Lie group acting properly
on
$M$. Let $X \in \mathfrak{X} (M) ^G$ and $F _t$ be the its
(necessarily equivariant) flow. For any isotropy subgroup $H$ of the
$G$-action on $M$, the
$H$-{\bfi  isotropy type submanifold} 
$M _H:=\{m \in M\mid G _m=H\}$ is preserved by the flow $F _t$. 
This property is known as the {\bfi law of
conservation of isotropy\/}. The
properness of the action guarantees that $G_m $ is compact and
that the (connected components of) $M _H$ are  embedded submanifolds
of $M$ for any closed subgroup  $H$ of $G $. The manifolds $M_H $
are,  in general, not closed in $M $. Moreover, the quotient group
$N(H)/H $ (where $N(H)$ denotes the normalizer of $H$ in
$G$) acts freely and properly on $M _H$. Hence, if
$\pi_H: M _H
\rightarrow M _H/(N(H)/H) $ denotes the projection onto orbit space
and $i _H: M _H \hookrightarrow M $ is the injection, the vector
field
$X$ induces a unique vector field $X ^H $ on the quotient $M
_H/(N(H)/H)$ defined by 
$X ^H \circ \pi_H= T \pi _H \circ X \circ i _H$,
whose flow $F _t ^H $ is given by $F ^H_t\circ \pi _H= \pi_H
\circ F _t
\circ i _H $. We will refer to $X ^H \in \mathfrak{X} \left(M
_H/(N(H)/H)\right)$ as the $H$-{\bfi  isotropy type reduced vector
field} induced by $X$. 

This reduction technique has been widely exploited in handling specific
dynamical systems. When the symmetry group $G$ is compact and we are dealing with
a linear action the construction of the quotient $M _H/(N(H)/H) $
can be implemented in a very explicit and convenient manner by using
the invariant polynomials of the action and the theorems of Hilbert
and Schwarz-Mather. 

\section{Symplectic reduction}

{\bfi  Symplectic\/} or {\bfi Marsden-Weinstein reduction} is a
procedure that implements symmetry reduction for the symmetric
Hamiltonian systems defined on a symplectic manifold $(M, \omega)$. The
particular case in which the symplectic manifold is a cotangent bundle is
dealt with separately in~\cite{cotangent bundle reduction}.  We recall
that the {\bfi  Hamiltonian vector field} $X _h
\in 
\mathfrak{X} (M) $ associated to the  {\bfi  Hamiltonian function} $h \in  C^\infty(M)$
is uniquely determined by the equality $\omega(X _h, \cdot )= \mathbf{d}
h$. In this context, the symmetries $\Phi: G \times  M \rightarrow M $ of
interest are given by {\bfi  symplectic} or {\bfi  canonical
transformations}, that is, $\Phi_g ^\ast  \omega= \omega $, for any $g
\in G $. For canonical actions each $G$-invariant function $h \in
C^\infty(M)^{G} $ has an associated $G$-symmetric  Hamiltonian vector
field $X _h$. A Lie algebra action $\varphi$ is called {\bfi  symplectic} 
or {\bfi  canonical}, if $\boldsymbol{\pounds}_{ \phi( \xi)} \omega = 0 $
for all $\xi\in \mathfrak{g}$, where $\boldsymbol{\pounds} $ denotes the
Lie derivative  operator. If the Lie algebra action is induced from a
canonical Lie group action by taking its infinitesimal generators, then
it is also canonical.

\medskip

\noindent {\bf Momentum maps.} The symmetry reduction described in the previous section
for general vector fields does not produce a well-adapted answer for symplectic
manifolds $(M, \omega)$ in the sense that the reduced spaces $M _H/ (N(H)/H) $
are, in general, not symplectic. To solve this problem one has to use
the conservation laws associated to the canonical action, that often
appear as {\bfi  momentum maps}.

Let $G$ be a Lie group acting canonically on the symplectic 
manifold $(M,\,\omega)$. Suppose that for any $\xi\in\mathfrak{g}$,
the vector field $\xi_M$ is  Hamiltonian, with
Hamiltonian function $\mathbf{J}^\xi\in C^{\infty}(M)$ and that  $\xi \in
\mathfrak{g}
\mapsto \mathbf{J}^\xi \in  C^\infty(M) $  is linear. The map
$\mathbf{J}:M\rightarrow\mathfrak{g}^\ast$ defined by the relation
$\langle\mathbf{J}(z),\,\xi\rangle=\mathbf{J}^\xi(z)$,
for all $\xi\in\mathfrak{g}$ and $z\in M$, is called a 
{\bfi momentum map\/}  of the $G$-action. 
Momentum maps, if they exist, are  determined up to a constant in
$\mathfrak{g}^\ast$ for any connected component of $M $. 

\begin{examples}
\normalfont
\label{Examples: linear momentum.}
\noindent\textbf{(i) Linear momentum.}
The phase space of an $N$-particle system is the cotangent space $T^\ast 
\mathbb{R}^{3 N}$ endowed with its canonical symplectic structure. The
additive group $\mathbb{R}^3$, whose Lie algebra is Abelian and is also
equal to $\mathbb{R}^3$, acts canonically on it by spatial translation on
each factor: $\boldsymbol{v}\cdot
(\boldsymbol{q}_i,\,\boldsymbol{p}^i)=(\boldsymbol{q}_i+
\boldsymbol{v},\,\boldsymbol{p}^i)$,
with $i=1,\ldots,\,N$. This action has an
associated momentum map $\mathbf{J}:T^\ast  \mathbb{R}^{3N}
\rightarrow \mathbb{R}^3$, where we identified the dual of $\mathbb{R}^3$
with itself using the Euclidean inner product, that coincides with the
classical {\bfi linear momentum\/}
$\mathbf{J}(\boldsymbol{q}_i,\,\boldsymbol{p}^i) =
\sum_{i=1}^{N}\boldsymbol{p}_i$.

\medskip

\noindent \textbf{(ii) Angular momentum.}
Let ${\rm SO}(3)$ act on $\mathbb{R}^3$ and then, by lift, on
$T^\ast \mathbb{R}^3$, that is,
$A\cdot(\boldsymbol{q},\,\boldsymbol{p})=(A\boldsymbol{q},\,
A\boldsymbol{p})$.
This action is canonical and has as associated momentum map
$\mathbf{J}: T^\ast  \mathbb{R}^{3} \rightarrow \mathfrak{so}(3)^\ast
\cong \mathbb{R}^3$ the classical {\bfi angular momentum\/}
$\mathbf{J}(\boldsymbol{q},\,\boldsymbol{p}) =
\boldsymbol{q}\times\boldsymbol{p}$.

\medskip

\noindent \textbf{(iii) Lifted actions on cotangent bundles.} 
The previous two examples are particular cases of the following situation. 
Let $\Phi: G \times M \rightarrow M $ be a smooth Lie group action. The
(left) {\bfi cotangent lifted action\/} of $G$ on $T^\ast Q$ is given
by $g \cdot \alpha_q : = T ^\ast_{g \cdot q} \Phi_{g ^{-1}}( \alpha_q) $
for $g \in G $ and $\alpha_q \in T ^\ast Q $. Cotangent lifted
actions preserve the canonical one-form on $T ^\ast Q $ and hence are
canonical. They admit an associated momentum map $\mathbf{J}:T^\ast 
Q\rightarrow\mathfrak{g}^\ast$ given by
$\langle\mathbf{J}(\alpha_q), \xi\rangle = \alpha_q(\xi_Q(q))$,
for any $\alpha_q\in T^\ast  Q$ and any $\xi\in \mathfrak{g}$. 

\medskip

\noindent\textbf{(iv) Symplectic linear actions.}
Let $(V,\,\omega)$ be a symplectic linear space and let $G$ be 
a subgroup of the linear symplectic group, acting naturally on
$V$. By the choice of $G$ this action
is canonical and has a momentum map given by 
$\langle\mathbf{J}(v),\,\xi\rangle=\frac{1}{2}\omega(\xi_V(v),\,v)$,
for $\xi\in \mathfrak{g}$ and $v\in V$  arbitrary.
\end{examples}

\noindent {\bf Properties of the momentum map.} The main feature of
the momentum map that makes it of interest for use in reduction is that it encodes
conservation laws for $G$-symmetric Hamiltonian systems. {\bfi  Noether's 
Theorem} states that the momentum map is a constant of the motion for the
Hamiltonian vector field $X _h$ associated to any $G$-invariant function
$h\in C^{\infty}(M)^{G}$.

The derivative $T \mathbf{J} $ of the momentum map satisfies the
following two properties: ${\rm range }\,(T _m \mathbf{J})=
(\mathfrak{g}_m) ^{\circ}$ and $\ker T _m
\mathbf{J}=(\mathfrak{g}\cdot m )^{\omega} $, for any $m \in  M $, 
where  $(\mathfrak{g}_m) ^{\circ} $ denotes the annihilator in
$\mathfrak{g}^\ast$ of the isotropy subalgebra $\mathfrak{g}_m $ of $m $,
$ \mathfrak{g}\cdot m:= T _m (G \cdot m) = \{ \xi_M(m) \mid \xi\in
\mathfrak{g}\}$ is the tangent space at
$m$ to the $G$-orbit that contains this point, and $(\mathfrak{g} \cdot m)
^\omega $ is the symplectic orthogonal space to $\mathfrak{g}\cdot m $ in
the symplectic vector space $( T_m M, \omega(m) ) $. The first relation is
sometimes called the {\bfi  bifurcation lemma\/}
since it establishes a link between the symmetry of a point and the rank
of the momentum map at that point.
 
The existence of the momentum map for a given canonical action is not guaranteed.
A momentum map exists if and only if the linear map $ \rho: [\xi] \in 
\mathfrak{g}/[\mathfrak{g}, \mathfrak{g}] \mapsto [ \omega( \xi_M, \cdot
) ] \in  H ^1 (M, \mathbb{R})$ is identically zero. Thus if $H ^1 (M,
\mathbb{R})=0 $ or $\mathfrak{g}/[ \mathfrak{g}, \mathfrak{g}] = H ^1
(\mathfrak{g}, \mathbb{R})=0 $ then $\rho\equiv 0 $. In
particular, if $\mathfrak{g} $ is semisimple, the First Whitehead Lemma
states that $H ^1 (\mathfrak{g}, \mathbb{R})=0 $ and therefore a 
momentum map always exists for canonical semisimple Lie algebra
actions.

A natural
question to ask is when the map
$(\mathfrak{g},[\cdot , \cdot ])
\rightarrow (C^\infty(M), \{ \cdot , \cdot \})$ defined by $ \xi \mapsto \mathbf{J}^\xi
$,
$\xi\in
\mathfrak{g}$, is a Lie algebra homomorphism, that is,
$
\mathbf{J}^{[\xi,\,\eta]}=\{\mathbf{J}^\xi,\,\mathbf{J}^\eta\}$, $\xi,
\eta \in \mathfrak{g}$. Here $\{ \cdot , \cdot \}: C^\infty(M)\times
C^\infty(M) \rightarrow C^\infty(M) $ denotes the Poisson bracket
associated to the symplectic form $\omega $ of
$M$ defined by $\{f,h\}:= \omega(X _f, X _h)$, $f,h \in  C^\infty(M) $.
This is the case if and only if
$T_z\mathbf{J}(\xi_M(z))=-\operatorname{ad}^\ast _\xi\mathbf{J}(z)$,
for any $\xi\in \mathfrak{g} $, $z \in M $,  where  $\mbox{\rm
ad}^\ast $ is the dual of the adjoint representation $\mbox{\rm ad}:
( \xi, \eta) \in \mathfrak{g} \times  \mathfrak{g} \mapsto [ \xi, \eta]
\in  \mathfrak{g} $  of $\mathfrak{g}$ on itself.
A momentum map that satisfies this relation in called {\bfi
infinitesimally equivariant\/}. 
The reason behind this terminology is that this is the infinitesimal
version of {\bfi  global} or {\bfi  coadjoint equivariance}:
 ${\bf J}$ is $G$-{\bfi  equivariant\/} if
$\operatorname{Ad}^\ast _{g^{-1}}\circ\mathbf{J}=\mathbf{J}\circ\Phi_g$
 or, equivalently,
$\mathbf{J}^{\operatorname{Ad}_g\xi}(g\cdot z)=\mathbf{J}^\xi(z)$,
for all $g\in G$, $\xi\in\mathfrak{g}$, and $z\in M$; $\mbox{\rm
Ad}^\ast$  denotes the dual of the adjoint representation $ \mbox{\rm
Ad}$ of $G$ on $\mathfrak{g}$. Actions admitting infinitesimally
equivariant momentum maps are called {\bfi Hamiltonian actions\/}
\index{action!Hamiltonian}%
\index{Hamiltonian!action}%
and Lie group actions 
with coadjoint equivariant momentum maps are called
{\bfi globally Hamiltonian actions\/}.
\index{action!globally Hamiltonian}%
\index{globally Hamiltonian action}%
If the symmetry group $G$ is connected then global and infinitesimal
equivariance of the momentum map are equivalent concepts.  If $ \mathfrak{g}$ acts
canonically on $(M, \omega )$ and $H ^1 (\mathfrak{g}, \mathbb{R})=\{0\} $ then this
action admits at most one infinitesimally equivariant momentum map.

Since momentum maps are not uniquely defined, one may ask
whether one can choose them to be equivariant. It turns out that if the momentum map 
is associated to the action of a compact Lie group  this can
always be done.  Momentum maps of cotangent lifted actions are also
equivariant as are momentum maps defined  by symplectic
linear actions. Canonical actions of semisimple  Lie algebras  on 
symplectic manifolds  admit infinitesimally equivariant momentum maps,
since the Second Whitehead Lemma  states that $H^2(\mathfrak{g},
\mathbb{R}) = 0$ if $\mathfrak{g}$ is semisimple. We shall identify below
a specific element of $H^2(\mathfrak{g}, \mathbb{R})$ which is the
obstruction to the equivariance of a momentum map (assuming it exists).

Even though, in general, it is not possible
to choose a coadjoint equivariant momentum map, it turns out that when the symplectic
manifold is connected there is an affine action on the dual of the Lie
algebra with respect to which the momentum map is equivariant. 
Define the {\bfi non-equivariance one-cocycle\/}
associated to $\mathbf{J}$ as the map
$\sigma:G\longrightarrow\mathfrak{g}^\ast$ given by
$g\longmapsto\mathbf{J}(\Phi_g(z))-\operatorname{Ad}^\ast _{g^{-1}}(\mathbf{J}(z))$.
The connectivity of $M$ implies that the right hand side  of this 
equality is independent of the point $z \in M $. In addition,  $\sigma$
is a (left) $\mathfrak{g}^\ast$-valued one-cocycle  on $G$ with respect to
the coadjoint representation of $G$ on
$\mathfrak{g}^\ast$, that is, $\sigma(gh) = \sigma(g) +
\operatorname{Ad}^\ast_{g ^{-1}} \sigma(h) $ for all $g, h \in G $.
Relative to the {\bfi affine action\/} $\Theta:G\times
\mathfrak{g}^\ast\longrightarrow\mathfrak{g}^\ast $ given by
$(g,\,\mu)\longmapsto\operatorname{Ad}^\ast _{g^{-1}}\mu+\sigma(g)$,
the momentum map $\mathbf{J}$ is equivariant. The {\bfi  Reduction
Lemma}, the main technical ingredient in the  proof of the reduction
theorem, states that for any $m \in M $ we have
\begin{equation*}
\label{reduction lemma in two equalities symplectic}
\mathfrak{g}_{\mathbf{J}(m)} \cdot m = \mathfrak{g}\cdot m\cap \ker T _m
\mathbf{J}=\mathfrak{g}\cdot m\cap (\mathfrak{g}\cdot m)^\omega,
\end{equation*}
where $\mathfrak{g}_{\mathbf{J}(m)} $  is the Lie algebra of the isotropy
group $G _{\mathbf{J}(m)}$ of $\mathbf{J}(m)
\in  \mathfrak{g}^\ast$ with respect to the affine action of $G $ on $\mathfrak{g}^\ast$
induced by the non-equivariance one-cocycle of  $\mathbf{J} $.

\medskip

\noindent \textbf{The Symplectic Reduction Theorem.} The symplectic reduction procedure that we
now present consists of constructing a new symplectic manifold out of a given
symmetric one in which the  conservation laws encoded in the form of a momentum map and
the degeneracies associated to the symmetry have been eliminated. This strategy allows
the reduction of a symmetric Hamiltonian dynamical system to a  dimensionally smaller
one. This reduction procedure preserves the symplectic category, that
is, if we start with a Hamiltonian system on a symplectic manifold, the reduced system is
also a Hamiltonian system on a symplectic manifold. The reduced symplectic manifold
is usually referred to as the {\bfi  symplectic} or  {\bfi
Marsden-Weinstein reduced space}.

\begin{theorem}
\label{symplectic point reduction}
Let $\Phi:G \times M \rightarrow M $ be a free proper canonical action of the Lie 
group $G$ on the connected symplectic manifold $(M,\omega)$.
Suppose that this action has an associated momentum map
$\mathbf{J} :M\rightarrow\mathfrak{g}^\ast$, with non equivariance
one-cocycle
$\sigma:G \rightarrow \mathfrak{g}^\ast$. Let $\mu\in\mathfrak{g}^\ast$ be a  value of
$\mathbf{J}$ and denote by $G_\mu$ the isotropy of $\mu$ under the affine action of $G$
on $\mathfrak{g}^\ast$. Then:
\begin{description}
\item[(i)] The space  $M_\mu:=\mathbf{J}^{-1}(\mu)/G_\mu$ is a regular quotient manifold and,
moreover, it is a symplectic  manifold with symplectic form $\omega_\mu$ uniquely 
characterized by the relation
\begin{equation*}
\label{mwredform}
\pi_\mu ^\ast\omega_\mu =i_\mu ^\ast \omega.
\end{equation*}
The maps $i _\mu: \mathbf{J}^{-1}(\mu) \hookrightarrow M $ and $\pi_{\mu}:
\mathbf{J}^{-1}(\mu)  \rightarrow \mathbf{J}^{-1}(\mu) / G _\mu $ denote the inclusion
and the projection, respectively.  The pair $(M_\mu,\, \omega _\mu)$ 
is called
the {\bfi symplectic point reduced space}.
\index{point!reduced space}%
\index{reduced!space!point}%
\index{space!point reduced}%

\item[(ii)] Let $h\in C^\infty(M)^G$ be a $G$-invariant
Hamiltonian. The flow $F_t$ of the Hamiltonian vector field $X_h$ leaves 
the connected components of $\mathbf{J}^{-1}(\mu)$ invariant and
commutes with the $G$-action, so it induces a flow $F_t^\mu$
on $M_\mu$ defined by
$
\pi_\mu \circ F_t \circ i_\mu=F_t^\mu \circ \pi_\mu$.
\item[(iii)] The vector field generated by the flow $F_t^\mu$ on $(M_\mu,\,\omega_\mu)$
is Hamiltonian with associated {\bfi reduced Hamiltonian function\/}
\index{Hamiltonian!reduced}
\index{reduced!Hamiltonian}
$h_\mu\in C^\infty(M_\mu)$ defined by
$h_\mu \circ \pi_\mu=h \circ i_\mu$.
The vector fields $X_h$ and $X_{h_\mu}$ are $\pi_\mu$-related.
The triple $(M_\mu,\,\omega_\mu, h_\mu)$ is called the {\bfi reduced Hamiltonian
system\/}.
\index{Hamiltonian!system!reduced}
\index{reduced!Hamiltonian!system}

\item[(iv)] Let $k\in C^\infty(M)^G$ be another $G$-invariant 
function. Then $\{h,\,k\}$ is also $G$-invariant and
$\{h,\,k\}_\mu=\{h_\mu,\,k_\mu\}_{M_\mu}$, where $\{ \cdot ,\cdot\} _{M _\mu}$ denotes
the Poisson bracket associated to the symplectic form  $\omega_\mu $ on $M _\mu$.
\end{description}
\end{theorem}

\medskip

\noindent \textbf{Reconstruction of dynamics.} 
We pose now the question converse to the  reduction of a Hamiltonian system. 
Assume that an integral curve $c_\mu(t)$
of the reduced Hamiltonian system $X_{h_\mu}$ on $(M_\mu, \omega_\mu)$ is known. Let
$m_0
\in
\mathbf{J}^{-1}(\mu)$ be given. One can determine from this data the integral curve of
the Hamiltonian system $X_h$ with initial condition $m_0$. In other words,  one can 
{\bfi reconstruct\/} 
the solution of the given system knowing the corresponding reduced solution. 
The general method of reconstruction is the following. Pick a smooth curve $d(t)$ in 
$\mathbf{J}^{-1}(\mu)$ such that $d(0) = m_0$ and $\pi_\mu(d(t)) = c_\mu(t)$. Then, if
$c(t)$ denotes the integral curve of $X_h$ with $c(0) = m_0$, we can write $c(t) = g(t)
\cdot d(t)$ for some smooth curve $g(t)$ in $G_\mu$ that is obtained in
two steps. First, one finds  a smooth
curve
$\xi(t)$ in
$\mathfrak{g}_\mu$ such that 
$\xi(t)_M(d(t)) = X_h(d(t)) - \dot{d}(t)$. With the 
$\xi(t) \in \mathfrak{g}_\mu$ just obtained, one solves
the non-autonomous differential equation 
$\dot{g}(t) = T_e L_{g(t)} \xi(t)$ on $G_\mu$ with $g(0) = e$.
\medskip

\noindent \textbf{The orbit formulation of the Symplectic Reduction Theorem.} There is 
an alternative approach to the reduction theorem  which
consists of choosing as numerator of the symplectic reduced space the group invariant
saturation of the level sets of the momentum map. This option
produces as a result a space that is symplectomorphic to the
Marsden-Weinstein quotient but presents the  advantage of being more
appropriate in the context of quantization problems. Additionally, this
approach makes easier the comparison of the symplectic reduced spaces
corresponding to different values of the momentum map which is important
in the context of Poisson reduction~\cite{Poisson reduction}. In carrying
out this construction one needs to use the natural symplectic structures
that one can define on the orbits of the affine action of a group on the
dual of its Lie algebra and that we now quickly review.

Let  $G$ be a Lie group, $\sigma:G\rightarrow \mathfrak{g}^\ast$  a
coadjoint one-cocycle, and  $\mu \in \mathfrak{g}^\ast$.  Let 
$\mathcal{O}_\mu$ be the orbit through $\mu$ of the affine
$G$-action on $\mathfrak{g}^\ast$ associated to $\sigma$. If  $\Sigma:
\mathfrak{g} \times \mathfrak{g} \rightarrow  \mathbb{R} $ defined by
$\Sigma(\xi, \eta) : = \left.\frac{d}{dt}\right|_{t=0}
\left\langle\sigma( \exp( t \xi) , \eta \right\rangle$ is a real valued 
Lie algebra two-cocycle (which is always the case if $\sigma$ is the
derivative of a smooth real valued group two-cocycle or if $\sigma$ is
the non-equivariance one-cocycle of a momentum map), that is, $\Sigma:
\mathfrak{g}\times \mathfrak{g}\rightarrow \mathbb{R}$ is skew symmetric
and $\Sigma([ \xi, \eta], \zeta) + \Sigma([ \eta, \zeta], \xi) +
\Sigma([
\zeta, \xi], \eta) = 0 $ for all $\xi, \eta, \zeta\in \mathfrak{g}$, 
then the affine orbit
$\mathcal{O}_\mu$ is a symplectic manifold with $G$-invariant symplectic
structure
$\omega_{\mathcal{O}_\mu}^\pm$ given by
\begin{equation}
\label{kks symplectic form explicit}
\omega_{\mathcal{O}_\mu}^\pm(\nu)(\xi_{\mathfrak{g}^\ast}(\nu),\,
\eta_{\mathfrak{g}^\ast}(\nu))
=\pm\langle\nu,\,[\xi,\,\eta]\rangle\mp\Sigma(\xi,\,\eta),
\end{equation}
for arbitrary $\nu\in \mathcal{O}_\mu$, and $\xi,\,\eta\in\mathfrak{g}$. 
The symbol
$\xi_{\mathfrak{g}^\ast}(\nu):= -{\rm ad}_\xi^\ast
\nu+\Sigma(\xi,\,\cdot)$
denotes the infinitesimal generator of the affine action on
$\mathfrak{g}^\ast$ associated to
$\xi\in \mathfrak{g}$. The
symplectic structures $\omega_{\mathcal{O}_\mu}^\pm$ on $\mathcal{O}_\mu$
are  called the $\pm$--{\bfi orbit} or {\bfi Kostant-Kirillov-Souriau
(KKS) symplectic forms}. 
 
This symplectic form can be obtained from Theorem~\ref{symplectic point reduction} by
considering the symplectic reduction of the cotangent bundle
$T ^\ast  G $ endowed with the {\bfi  magnetic symplectic structure}
$\overline{\omega _\Sigma}:=
\omega _{{\rm can}}- \pi ^\ast  B _\Sigma  $, where $\omega _{{\rm can}}$  is the
canonical symplectic form  on $T ^\ast G $,
$\pi :T ^\ast  G \rightarrow  G $  is the projection onto the base, and $B _\Sigma
\in  \Omega ^2(G)^G $ is a left invariant  two-form on $G$  whose value at the
identity is the Lie algebra two-cocycle $\Sigma: \mathfrak{g} \times   \mathfrak{g}
\rightarrow \mathbb{R}$. Since $\Sigma $  is a cocycle it follows that $B _\Sigma 
$ is closed and hence $\overline{\omega _\Sigma}$ is a symplectic  form. Moreover, the
lifting of the left translations on $G$ provides a canonical $G$-action on $T ^\ast G $
that has a momentum map given by $\mathbf{J}(g, \mu)= \Theta(g, \mu) $, $(g, \mu) \in 
G \times  \mathfrak{g}^\ast \simeq T ^\ast G $, where the trivialization $G \times 
\mathfrak{g}^\ast \simeq T ^\ast G $ is obtained via left translations. Symplectic
reduction using these ingredients yields symplectic reduced spaces that are naturally
symplectically diffeomorphic to the affine orbits $\mathcal{O}_{\mu}$ with the
symplectic form~(\ref{kks symplectic form explicit}).

\begin{theorem}[Symplectic orbit reduction]
\label{orbitreductionregular}
Let $\Phi:G \times M \rightarrow M $ be a free proper canonical action of the Lie 
group $G$ on the connected symplectic manifold $(M,\omega)$.
Suppose that this action has an associated momentum map
$\mathbf{J}:M\rightarrow\mathfrak{g}^\ast$, with non equivariance
one-cocycle $\sigma:G \rightarrow
\mathfrak{g}^\ast$. Let $\mathcal{O}_\mu:=G \cdot \mu\subset\mathfrak{g}^\ast$ be the
$G$-orbit of the point
$\mu \in \mathfrak{g}^\ast$ with respect to the affine action of $G$ on
$\mathfrak{g}^\ast$ associated to $\sigma$. Then
the set $M_{\mathcal{O}_\mu}:=\mathbf{J}^{-1}(\mathcal{O}_\mu)/G$ is a regular quotient
symplectic manifold with the symplectic form $\omega_{ \mathcal{O}_{\mu}} $ uniquely
characterized by the relation
$
i_{\mathcal{O}_\mu}^\ast\omega=\pi_{\mathcal{O}_\mu}^\ast\omega_{\mathcal{O}_\mu}
+\mathbf{J}^\ast_{\mathcal{O}_\mu}\omega_{\mathcal{O}_\mu}^+$,
where $\mathbf{J}_{\mathcal{O}_\mu}$ is the restriction of $\mathbf{J}$ to $\mathbf{J}^{-1}(\mathcal{O}_\mu)$ and
$\omega_{\mathcal{O}_\mu}^+$ is the $+$--symplectic structure on the affine
orbit $\mathcal{O}_{\mu} $. The maps
$ i _{\mathcal{O}_{\mu}}:
\mathbf{J}^{-1}(\mathcal{O}_{\mu}) \hookrightarrow
M $ and $\pi_{\mathcal{O}_\mu}:\mathbf{J}^{-1}(\mathcal{O}_\mu)\rightarrow M_{\mathcal{O}_\mu}$ are natural injection and  the
projection, respectively. The pair $(M_{\mathcal{O}_\mu},\, \omega_{\mathcal{O}_\mu})$
is called the {\bfi symplectic orbit reduced space}.
Statements similar to {\bf (ii)} through {\bf (iv)} in Theorem~\ref{symplectic point
reduction} can be formulated for the orbit reduced spaces  $(M_{\mathcal{O}_\mu},\,
\omega_{\mathcal{O}_\mu})$. 
\end{theorem}

We emphasize that given a momentum value $\mu \in \mathfrak{g}^\ast$, the reduced spaces
$M _\mu
$ and
$M _{\mathcal{O}_{\mu}}$ are symplectically diffeomorphic via the projection to the
quotients of the inclusion $\mathbf{J}^{-1}(\mu)\hookrightarrow
\mathbf{J}^{-1}(\mathcal{O}_{\mu}) $.

Reduction at a general point can be replaced by reduction at zero at the expense of
enlarging the manifold by the affine orbit.
Consider the canonical diagonal action of $G$ on the symplectic difference $M \ominus
\mathcal{O}_\mu^+$ which is the manifold $M \times  \mathcal{O}_{\mu} $ with 
the symplectic form $\pi _1 ^\ast  \omega - \pi _2 ^\ast \omega _{\mathcal{O}_{\mu}^+}
$, where $\pi_1: M
\times
\mathcal{O}_\mu \rightarrow M$ and $\pi_2: M \times \mathcal{O}_\mu \rightarrow 
\mathcal{O}_\mu $ are the projections. A momentum map for this action is given by
$\mathbf{J}
\circ
\pi_1 -
\pi_2 :M \ominus \mathcal{O}_\mu^+ \rightarrow \mathfrak{g}^\ast$. Let $(M\ominus \mathcal{O}_\mu^+)_0 : =
((\mathbf{J} \circ \pi_1 - \pi_2)^{-1}(0)/G, (\omega \ominus
\omega_{\mathcal{O}_\mu}^+)_0)$ be the symplectic point reduced space at zero.

\begin{theorem}[Shifting theorem]
\label{shifting theorem free}
\index{shifting theorem}
\index{theorem!shifting}
Under the hypotheses of the Symplectic Orbit Reduction
Theorem~{\rm \ref{orbitreductionregular}}, the symplectic orbit reduced space
$M_{\mathcal{O}_\mu}$,  the point reduced spaces
$M_\mu$, and $(M\ominus \mathcal{O}_\mu^+)_0$ are symplectically diffeomorphic.
\end{theorem}

\section{Singular reduction}

In the previous section we carried out symplectic reduction for free and proper actions.
The freeness guarantees via the bifurcation lemma that the momentum map ${\bf J}$ is a
submersion and hence the level sets $\mathbf{J}^{-1}(\mu) $ are smooth manifolds.
Freeness and properness ensures that the orbit spaces $M _\mu:=
\mathbf{J}^{-1}(\mu)/ G_{\mu}$ are regular quotient manifolds. The theory of singular
reduction studies the properties of the orbit space $M _\mu $ when the hypothesis
on the freeness of the action is dropped.  The main
result in this situation shows that these quotients are  
symplectic Whitney stratified spaces in the sense that the strata are symplectic
manifolds in a very natural way; moreover, the local properties of this
Whitney stratification make it into what is called a {\bfi  cone space}.  This statement
is referred to as the {\bfi  Symplectic Stratification Theorem} and adapts to the
symplectic symmetric context the stratification theorem of the orbit space of a proper
Lie group action by using its orbit type manifolds.
In order to present this  result we review the necessary definitions
and results on stratified spaces.

\medskip

\noindent {\bf Stratified spaces.} Let ${\cal Z}$ be a locally finite partition of the
topological space $P$ into smooth manifolds $S _i \subset P $, $i \in
I$. We assume that the manifolds $S _i\subset P$,
$i\in I$, with their manifold topology are locally closed topological
subspaces  of $P$. The pair $(P, {\cal Z})$ is a {\bfi  decomposition}
\index{decomposition}%
of $P$ with {\bfi  pieces}
\index{piece!of a decomposition}%
\index{decomposition!pieces}%
in ${\cal Z}$ when the following condition is satisfied:
\begin{description}
\item [(DS)] If $R, S \in {\cal Z}$ are such that $R \cap \bar{S} \neq
\varnothing$, then
$R \subset \bar{S} $. In this case we write $R \preceq S $. If, in addition, $R \neq S$
we say that
$R$ is {\bfi  incident } 
\index{incident}%
to $S$ or that it is a {\bfi  boundary piece} 
\index{boundary!piece}%
\index{piece!boundary}%
of $S$ and write $R \prec S$.  
\end{description}
Condition {\bf (DS) } is called the {\bfi  frontier condition}
and the pair $(P, {\cal
Z})$ is called a {\bfi  decomposed space}.
The {\bfi  dimension} 
\index{dimension!of a decomposed space}%
\index{decomposed!space!dimension}%
of $P$ is defined as $
\dim P=\sup \{\dim S _i\mid S_i \in {\cal Z}\}$. If $k \in \mathbb{N}$, the
$k$-{\bfi  skeleton} 
\index{skeleton of a decomposed space}%
\index{decomposed!space!skeleton}%
$P ^k $ of $P$ is the union of all the pieces of dimension smaller than or equal to
$k$; its topology is the relative topology induced by $P$. The
{\bfi  depth}
\index{depth}%
${\rm d p} (z) $ of any $z \in (P, {\cal Z})$ is defined as 
\[
{\rm d p}(z):=\sup \{k \in \mathbb{N}\mid \exists\  S _0, S _1, \ldots, S _k \in {\cal
Z}\ {\rm with }\ z \in S _0 \prec S _1 \prec \ldots \prec S _k\}.
\]
Since for any two elements $x,y \in \mathcal{S} $ in the same piece
$\mathcal{S} \in P $  we have ${\rm d p}(x)= {\rm d p} (y) $, the
depth ${\rm d p} (S) $ of the piece
$S $ is well defined by ${\rm d p} (S):= {\rm d p} (x) $, $x \in \mathcal{S} $. Finally,
the depth ${\rm d p} (P) $ of $(P, {\cal Z} ) $ is defined by ${\rm d p} (P) :=\mbox{\rm
sup}\{ {\rm d p} (S)\mid S \in {\cal Z}\}$.

A
continuous mapping
$f : P
\rightarrow Q
$ between the decomposed spaces
$(P, {\cal Z})$ and $(Q, {\cal Y}) $ is a {\bfi  morphism of decomposed spaces}
\index{morphism!of decomposed spaces}%
\index{decomposed!space!morphism}%
if for every piece $S \in {\cal Z} $, there is a piece $T \in {\cal Y} $ such that $f
(S)\subset  T $ and the restriction $f| _S: S \rightarrow T $ is smooth. If $(P, {\cal
Z})$ and $(P, {\cal T})$ are two decompositions of the same topological
space we say that
${\cal Z}$ is {\bfi  coarser}
\index{coarser decomposition}%
\index{decomposition!coarser}%
than ${\cal T}$ or that ${\cal T}$ is {\bfi  finer} 
\index{finer decomposition}%
\index{decomposition!finer}%
than ${\cal Z}$ if the identity mapping $(P, {\cal T})\rightarrow (P, {\cal Z})$ is a
morphism of decomposed spaces. A topological subspace $Q \subset P $ is a {\bfi 
decomposed subspace}
\index{decomposed!subspace}%
\index{subspace!decomposed}%
of $(P, {\cal Z}) $ if for all pieces $S \in {\cal Z} $, the intersection  $S\cap Q $ is
a submanifold of $S$ and the corresponding partition ${\cal Z} \cap Q $
forms a decomposition of $Q$.

Let  $P$ be a topological space and  $z \in P $. Two subsets
$A$ and  $B$ of $P$ are said to be {\bfi  equivalent}
\index{sets!equivalent}%
\index{equivalent!sets}%
at $z$ if there is an open neighborhood $U$ of $z$ such that $A\cap U= B\cap U $. This
relation constitutes an equivalence relation on the power set of $P$. The class of all
sets equivalent to a given subset $A$ at $z$ will be denoted by  $[A] _z$ and called the
{\bfi  set germ} 
\index{germ!set}%
\index{set!germ}%
of $A$ at $z$. If $A \subset B \subset P $ we say that $[A] _z$ is a {\bfi  subgerm}
\index{subgerm}%
of $[B] _z$, and denote $[A] _z \subset [B] _z $.

A {\bfi  stratification}
of the topological space $P$ is a map ${\cal S}$ that associates to any $z \in P $ the
set germ $\mathcal{S}(z) $ of a  closed subset of $P$ such that the following condition
is satisfied:
\begin{description}
\item [(ST)] For every $z \in P $ there is a neighborhood $U$ of $z$ and a
decomposition ${\cal Z}$ of $U$ such that for all
$y
\in  U $ the germ $ \mathcal{S} (y) $ coincides with the set germ of the piece of ${\cal
Z}$ that contains $y$.
\end{description}
The pair $(P , \mathcal{S}) $ is called a {\bfi  stratified space}.
\index{space!stratified}%
\index{stratified!space}%
Any decomposition of $P$ defines a stratification of  $P$ by associating to each of its
points the set germ of the piece in which it is contained. The converse
is, by definition, locally true. 

\medskip

\noindent {\bf The strata.} Two decompositions ${\cal Z}_1 $ and ${\cal Z}_2 $
of $P$ are said to be {\bfi equivalent}
\index{decomposition!equivalent}%
\index{equivalent!decomposition}%
if they induce the same stratification of $P$. If ${\cal Z}_1 $ and ${\cal Z}_2 $ are
equivalent decompositions of
$P$ then, for all $z \in P $, we have that ${\rm d p} _{{\cal Z}_1}(z)=
{\rm d p} _{{\cal Z}_2}(z) $. Any stratified space $(P, \mathcal{S})$ has
a unique  decomposition ${\cal Z}_{\mathcal{S}} $ associated with the
following maximality property: for any open subset $U \subset P $ and any
decomposition ${\cal Z}$ of $P$ inducing ${\cal S}$ over $U$, the
restriction of ${\cal Z}_{\mathcal{S}} $ to $U$ is coarser than the
restriction of ${\cal Z}$ to $U$. The decomposition ${\cal
Z}_{\mathcal{S}} $ is called the {\bfi  canonical decomposition}
\index{decomposition!associated to a stratification}%
\index{stratification!canonical associated decomposition}%
\index{canonical!decomposition associated to a stratification}%
associated to the stratification $(P, \mathcal{S})$. It is often denoted by ${\cal S}$
and its pieces are called the {\bfi  strata}
\index{stratum}%
of $P$. The local finiteness of the decomposition ${\cal Z}_{\mathcal{S}} $ implies that
for any stratum $S $ of $( P, \mathcal{S}) $ there are only finitely many strata $R$ with
$S \prec R $. In the sequel the symbol  $ \mathcal{S}$  in the stratification $( P,
\mathcal{S}) $ will denote both the map that associates to each point a set germ and the
set of pieces associated to the canonical decomposition induced by the stratification of
$P$.

\medskip

\noindent {\bf Stratified spaces with smooth structure.} 
Let  $(P , \mathcal{S}) $ be
a stratified space. A {\bfi  singular} or {\bfi  stratified chart}
\index{chart!singular}%
\index{chart!stratified}%
\index{singular!chart}%
\index{stratified!chart}%
of $P$ is a homeomorphism $\phi: U \rightarrow \phi(U) \subset \Bbb R^n $ from an open set
$U \subset P $ to a subset of $\Bbb R^n $ such that for every stratum
$S \in \mathcal{S}$ the image $\phi(U \cap S) $ is a submanifold of $\Bbb R^n$ and
the restriction $\phi|_{U\cap S}:U\cap S \rightarrow \phi(U\cap S) $ is a diffeomorphism.
Two singular charts $\phi: U \rightarrow \phi(U) \subset \Bbb R^n $ and $\varphi: V
\rightarrow \varphi(V) \subset \Bbb R^m $ are {\bfi  compatible}
\index{compatible!stratified charts}%
\index{chart!stratified!compatibility of}%
if for any $z \in U\cap V $ there exist an open neighborhood $W \subset U\cap V $ of 
$z$, a natural number $N \geq \max\{n, m\} $, open neighborhoods $O , O' \subset \Bbb R
^{N}$ of $\phi (U) \times \{0\} $ and $\varphi (V) \times \{0\} $, respectively, and a
diffeomorphism $\psi: O \rightarrow O' $ such that $i _m \circ \varphi| _W= \psi \circ  i
_n \circ \phi |_W$, where $i _n $ and $i _m$ denote the natural embeddings of $\Bbb R^n
$ and $\Bbb R^m$ into  $\Bbb R^N $ by using the first $n $ and $m$ coordinates,
respectively. The notion of {\bfi  singular} or {\bfi  stratified atlas}
\index{atlas!singular}%
\index{atlas!stratified}%
\index{stratified!atlas}%
\index{singular!atlas}%
is the natural generalization for stratifications of the concept of atlas existing for
smooth manifolds. Analogously, we can talk of compatible and maximal stratified atlases.
If the stratified space $(P, \mathcal{S}) $ has a well defined maximal atlas, then we say
that this atlas determines a {\bfi  smooth } or {\bfi  differentiable structure} 
\index{smooth!structure of a stratification}%
\index{differential structure of a stratification}%
on $P$. We will refer to $(P, \mathcal{S}) $ as a {\bfi  smooth stratified space}.

\medskip

\noindent {\bf The Whitney conditions.}  Let $M$ be a manifold and $R, S \subset M$  
two
submanifolds. We say that the pair $(R, S)$ satisfies the {\bfi  Whitney condition} {\bf
(A)}
at the point $z \in  R $
if the following condition is satisfied:
\begin{description}
\item [(A)] For any sequence of points $\{z _n\} _{n \in \Bbb N} $ in $S$ converging to
$z \in R $ for which the sequence of tangent spaces $\{T _{z _n} S\} _{n \in \Bbb N}
$ converges in the Grassmann bundle  of $\dim S $--dimensional subspaces of $TM $ to $
\tau\subset T _z M$, we have that $T _z R \subset \tau$.
\end{description}
Let $\phi: U \rightarrow \Bbb R^n$ be a smooth chart of $M$ around the point $z$. The
{\bfi  Whitney condition} {\bf
(B)}
\index{condition!Whitney (B)}%
\index{Whitney!condition (B)}%
\index{regular!(B)--pair}%
at the point $z \in  R $ with respect to the chart $(U,\phi)$ is given by the following
statement:
\begin{description}
\item [(B)] Let $\{x _n\} _{n \in \Bbb N} \subset R\cap U$ and $\{y _n\} _{n \in \Bbb N}
\subset S\cap U$ be two sequences with the same limit $z=\lim\limits _{n \rightarrow
\infty}x _n=\lim\limits _{n \rightarrow \infty}y
_n$ and such that
$x _n\neq y _n
$, for all
$n\in\Bbb N$. Suppose that the set of connecting lines $\overline{\phi(x _n)\phi(y
_n)}\subset \Bbb R^n$ converges in projective space to a line $L$ and that the sequence of
tangent spaces $\{T _{y _n} S\} _{n \in \Bbb N}
$ converges in the Grassmann bundle  of $\dim S $-dimensional subspaces
of $TM $ to $
\tau\subset T _z M$. Then, $(T _z \phi)^{-1}(L)\subset \tau$.
\end{description}
If the condition (A) (respectively (B)) is verified for every point $z \in R $, the
pair $(R, S) $ is said to satisfy the {\bfi  Whitney condition} \textbf{(A)} (respectively {\bf (B)}).
It can be verified that Whitney's condition (B) does not depend on the chart used to
formulate it. A
stratified space with smooth structure such that for every pair of strata Whitney's
condition (B) is satisfied is called a {\bfi  Whitney space}.

\medskip

\noindent {\bf Cone spaces and local triviality.} Let $P$ be a topological space. Consider
the equivalence relation $\sim $ in the product $P \times [0 , \infty) $ given by 
$(z,a)\sim (z', a') $ if and only if $a=a'=0 $.  We define the {\bfi  cone} 
\index{cone}%
$C P $ on $P$ as the quotient topological space $P \times [0 , \infty) / \sim $. If $P$
is a smooth manifold then the cone $C P $ is a decomposed space with two
pieces, namely,
$P
\times (0,
\infty)$ and the {\bfi  vertex} 
\index{vertex of a cone}%
\index{cone!vertex}%
which is
the class corresponding to any element of the form $(z,0)$,  $ z \in P $, that is, $P
\times \{0\} $. Analogously, if $(P , {\cal Z} ) $ is a decomposed (stratified) space then
the associated cone $C P $ is also a decomposed (stratified) space  whose pieces
(strata) are the vertex and the sets of the form $S \times (0, \infty)$, with $S \in {\cal
Z} $. This implies, in particular, that 
$
\dim C P=\dim P+ 1$ and 
${\rm d p}(CP)= {\rm d p} (P)+ 1$.

A stratified space $(P, {\cal S})$ is said to be {\bfi  locally trivial}
if for any $z \in P $ there exist a
neighborhood $U$ of $z$, a stratified space $(F, \mathcal{S}^F) $, a distinguished point
${\bf 0} \in F $, and an isomorphism of stratified spaces
\[\psi : U  \rightarrow (S\cap U)\times
F ,\] 
where $S$ is the stratum that contains $z$ and $\psi$ satisfies  $ \psi
^{-1}(y , {\bf 0})=y $, for all $y \in S\cap U $.
When $F $ is given by a cone $CL$
over a compact stratified space $L$ then $L$ is called the {\bfi  link}
\index{link}%
of $z$. 

An important corollary of
{\bfi  Thom's First Isotopy Lemma}
guarantees that every Whitney stratified space is locally trivial. A converse to this
implication needs the introduction of cone spaces. Their definition is
given by recursion on the depth of the space.
\begin{definition}
\label{definition cone space 13}
Let $m \in \mathbb{N} \cup \{ \infty, \omega\} $. A {\bfi  cone space}
\index{cone!space}%
\index{space!cone}%
of class $C ^m $ and depth  $0$ is the union of countably many $C ^m $  manifolds together
with the stratification whose strata are the unions of the connected components of equal
dimension. A cone space of class $C ^m $ and depth $d + 1 $, $d \in \mathbb{N} $, is a
stratified space $(P, \mathcal{S})$ with a $C ^m  $ differentiable  structure such that
for any $z \in P $ there exists a connected neighborhood $U $ of $z $, a compact cone
space $L$ of class $C ^m  $ and depth $d  $ called the {\bfi  link}, 
\index{link}%
and a stratified
isomorphism
\[
\psi:U \rightarrow (S\cap U) \times C L,
\]
where $S$ is the stratum that contains the point $z$, the map  $\psi$ satisfies 
$ \psi ^{-1}(y , {\bf 0})=y $, for all $y \in S\cap U $, and ${\bf 0 }$
is the vertex of the cone
$C L$.

If $m \neq 0 $ then $L$ is required to be embedded into a sphere via a fixed smooth
global singular chart $\varphi:L \rightarrow S ^l $  that determines the smooth structure
of $C L $. More specifically, the smooth structure of $C L $ is generated by the global
chart $\tau:[z,t] \in C L  \longmapsto t \varphi(z) \in 
\mathbb{R}^{l+1} $. The maps $\psi:U \rightarrow (S\cap U) \times C L$
and $\varphi:L \rightarrow S ^l $ are referred to as a {\bfi  cone chart}
\index{chart!cone}%
\index{cone!chart}%
and a {\bfi  link chart},
\index{chart!link}%
\index{link!chart}%
respectively. Moreover, if $m \neq 0 $ then $\psi $ and $\psi ^{-1} $ are required to  be
differentiable of class $C ^m  $  as maps between stratified spaces with a smooth
structure.
\end{definition}

\medskip

\noindent {\bf The Symplectic Stratification Theorem.} Let $(M, \omega)$ be a connected symplectic manifold acted canonically and properly upon by
a Lie group $G$. Suppose that this action has an associated momentum map
$\mathbf{J}: M \rightarrow \mathfrak{g}^\ast$ with non-equivariance
one--cocycle $\sigma:G
\rightarrow\mathfrak{g}^\ast$.  Let
$\mu\in\mathfrak{g}^\ast$ be a value of $\mathbf{J}$,  $G _\mu $  the isotropy subgroup
of 
$\mu$ with respect to the affine action $\Theta:G \times \mathfrak{g}^\ast \rightarrow
\mathfrak{g}^\ast$ determined by $\sigma$, and let $H\subset G$ be an isotropy subgroup
of the $G$-action on $M$. Let $M _H^z $ be the connected component of the
$H$-isotropy type manifold that contains a given element $z \in M $ such
that $\mathbf{J}(z)= \mu $ and let $G _\mu M _H^z$ be its $G
_\mu$-saturation. Then the following hold:
\begin{description}
\item[(i)] The set $\mathbf{J}^{-1}(\mu)\cap G _\mu M _H^z$ is a submanifold of $M$.
\item[(ii)] The set $M_{\mu}^{(H)}:=[\mathbf{J}^{-1}(\mu)\cap G _\mu M _H^z]/G _\mu$
has a unique quotient differentiable structure such that the canonical projection
$\pi_{\mu}^{(H)}:\mathbf{J}^{-1}(\mu)\cap G _\mu M _H^z\longrightarrow
M_{\mu}^{(H)}$ 
is a surjective submersion. 
\item[(iii)] There is a unique symplectic structure 
$\omega_{\mu}^{(H)}$ on $M_{\mu}^{(H)}$ characterized by
\[
i_{\mu}^{(H)\,*}\omega=\pi_{\mu}^{(H)\,*}\omega_{\mu}^{(H)},
\]
where $i_{\mu}^{(H)}:\mathbf{J}^{-1}(\mu )\cap G _\mu M _H^z\hookrightarrow M$ is
the natural inclusion. The pairs $(M _\mu^{ (H)}, \omega_\mu^{ (H)})$ will be
called {\bfi  singular symplectic point strata}.
\item[(iv)] Let $h\in C^{\infty}(M)^G$ be a $G$-invariant 
Hamiltonian. Then the flow $F_{t}$ of $X_{h}$ leaves the 
connected components of 
$\mathbf{J}^{-1}(\mu)\cap G _\mu M _H^z$ invariant 
and commutes with the $G_{\mu}$-action, so it induces a flow
$F_{t}^{\mu}$ on $M_{\mu}^{(H)}$ that is characterized by 
$
\pi_{\mu}^{(H)}\circ F_{t}\circ i_{\mu}^{(H)} =F_{t}^{\mu}\circ 
\pi_{\mu}^{(H)}$.
\item[(v)] The flow $F_{t}^{\mu}$ is Hamiltonian 
on $M_{\mu}^{(H)}$, with {\bfi reduced Hamiltonian}
function $h_{\mu}^{(H)}:M_{\mu}^{(H)}\rightarrow \mathbb{R}$ defined by
$
h_{\mu}^{(H)}\circ \pi_{\mu}^{(H)}=h\circ i_{\mu}^{(H)}$.
The vector fields $X_{h}$ and $X_{h_{\mu}^{(H)}}$ are 
$\pi_{\mu}^{(H)}$--related. 
\item[(vi)] Let $k:M\rightarrow \mathbb{R}$ be another 
$G$-invariant function. Then $\{h,k\}$ is also $G$-invariant and 
$
\{h,k\}_{\mu}^{(H)}=\{h_{\mu}^{(H)},k_{\mu}^{(H)}\}_{M_{\mu}^{(H)}}$,
where $\{\ ,\ \}_{M_{\mu}^{(H)}}$ denotes the 
Poisson bracket induced by the symplectic structure on $M_{\mu}^{(H)}$.
\end{description} 

\begin{theorem}[Symplectic Stratification Theorem]
The quotient $M _\mu:= \mathbf{J}^{-1}(\mu)/ G_{\mu} $ is a cone space  when considered
as a stratified space with strata $M _\mu^{(H)} $.
\end{theorem}

As was the case for regular reduction, this theorem can be also
formulated from the orbit reduction point of view. Using that
approach one can conclude that the orbit reduced spaces $M
_{\mathcal{O}_{\mu}} $ are cone spaces symplectically stratified by the
manifolds
$M_{\mathcal{O}_{\mu}}^{(H)}:=G \cdot (\mathbf{J}^{-1}(\mu )\cap M
_H^z)/G$ that have  symplectic structure uniquely determined by the
expression
\begin{equation*}
\label{expression orbit singular symplectic}
i_{\mathcal{O}_{\mu}}^{(H)\,*}\omega=\pi_{\mathcal{O}_{\mu}}^{(H)\,*}
\omega_{\mathcal{O}_{\mu}}^{(H)}+
\mathbf{J}_{\mathcal{O}_{\mu}}^{(H)\,\ast}\omega_{\mathcal{O}_{\mu}} ^+, 
\end{equation*}
where
$i_{\mathcal{O}_{\mu}}^{(H)}:G
\cdot (\mathbf{J}^{-1}(\mu )\cap M _H^z)\hookrightarrow M$ is the inclusion,
$\mathbf{J}_{\mathcal{O}_{\mu}}^{(H)}:G \cdot (\mathbf{J}^{-1}(\mu )\cap M _H^z)
\rightarrow \mathcal{O}_{\mu} $ is obtained by restriction of the momentum map
$\mathbf{J}$, and $\omega_{\mathcal{O}_{\mu}} ^+ $ is the $+ $--symplectic form on
$\mathcal{O}_{\mu}$. Analogous statements to \textbf{(i) - (vi)} above
with obvious modifications are valid.

\end{document}